\documentclass[12pt]{article}
\usepackage{amssymb,amsmath,bm}
\usepackage{mathrsfs}
\usepackage{constants}
\begin{document}
\newcommand{\bea}{\begin{eqnarray}}
\newcommand{\ena}{\end{eqnarray}}
\newcommand{\beas}{\begin{eqnarray*}}
\newcommand{\enas}{\end{eqnarray*}}
\newcommand{\beq}{\begin{equation}}
\newcommand{\enq}{\end{equation}}
\def\qed{\hfill \mbox{\rule{0.5em}{0.5em}}}
\newcommand{\bbox}{\hfill $\Box$}
\newcommand{\ignore}[1]{}
\newcommand{\ignorex}[1]{#1}
\newcommand{\wtilde}[1]{\widetilde{#1}}
\newcommand{\qmq}[1]{\quad\mbox{#1}\quad}
\newcommand{\qm}[1]{\quad\mbox{#1}}
\newcommand{\nn}{\nonumber}
\newcommand{\Bvert}{\left\vert\vphantom{\frac{1}{1}}\right.}
\newcommand{\To}{\rightarrow}

\newtheorem{theorem}{Theorem}[section]
\newtheorem{corollary}{Corollary}[section]
\newtheorem{conjecture}{Conjecture}[section]
\newtheorem{proposition}{Proposition}[section]
\newtheorem{lemma}{Lemma}[section]
\newtheorem{definition}{Definition}[section]
\newtheorem{example}{Example}[section]
\newtheorem{remark}{Remark}[section]
\newtheorem{case}{Case}[section]
\newtheorem{condition}{Condition}[section]
\newcommand{\pf}{\noindent {\it Proof:} }
\newcommand{\proof}{\noindent {\it Proof:} }

\title{Clubbed Binomial Approximation for the Lightbulb Process}

\author{Larry Goldstein\thanks{Department of Mathematics, University of Southern California} and Aihua Xia\thanks{Department of Mathematics and Statistics, The University of Melbourne}}
\maketitle

\abstract{In the so called lightbulb process, on days $r=1,\ldots,n$, out of $n$ lightbulbs, all initially off, exactly
$r$ bulbs selected uniformly and independent of the past have their status changed from off to on, or vice versa. With $W_n$ the number of bulbs on at the terminal time $n$ and $C_n$ a suitable clubbed binomial distribution,
$$
d_{{\mbox{\scriptsize \rm TV}}}(W_n,C_n) \le 2.7314 \sqrt{n} e^{-(n+1)/3} \quad \mbox{for all $n \ge 1$.}
$$
The result is shown using Stein's method.}

\section{Introduction}

The lightbulb process introduced by \cite{rrz} was motivated by a pharmaceutical study of the effect
of dermal patches designed to activate targeted receptors.
An active receptor will become inactive, and an inactive one active,
if it receives a dose of medicine released from the dermal patch. On each of $n$ successive
days $r=1,\ldots,n$ of the study, exactly $r$ randomly selected receptors will each receive one dose
of medicine from the patch, thus changing, or toggling, their status between the active and inactive states.
We adopt the more colorful language of ~\cite{rrz}, where receptors are represented by lightbulbs that
are being toggled between their on and off states.

Some fundamental properties of $W_n$, the number of light bulbs on at the end of day $n$, were derived in \cite{rrz}. For instance, Proposition~2 of \cite{rrz} shows that when $n(n+1)/2 = 0 \,\mbox{mod}\,2$, or, equivalently, when $n\,\mbox{mod}\,4 \in \{0,3\}$, the support of $W_n$ is a subset of even integers up to $n$, and that otherwise the support of $W_n$ is a set of odd integers up to $n$.
Further, in \cite{rrz}, the mean and variance of $W_n$ were computed, and based on numerical computations, an approximation of the distribution of $W_n$ by the `clubbed' binomial distribution was suggested.

To describe the clubbed binomial, let $Z_n$ be a binomial $\mbox{Bin}(n-1,1/2)$ random variable, and for $i \in \mathbb{Z}$ let $\pi_i^*=P(Z_n=i)$, that is
\beas
\pi_i^* = \left\{
\begin{array}{cc}
{n-1 \choose i} \left(\frac{1}{2}\right)^{n-1} &  \mbox{for $i=0,1,\ldots,n-1$,}\\
0 & \mbox{otherwise.}
\end{array}
\right.
\enas
Let $L_{1,n}$ and $L_{0,n}$ denote the set of all odd and even numbers in $\{0,1,\ldots,n\}$, respectively.
Define, for $m=0,1$,
$$\pi^m_i=\left\{\begin{array}{ll}
\pi_{i-1}^*+\pi_i^*,&\ i\in L_{m,n},\\
0,&\ i\not\in L_{m,n}.
\end{array}\right.
$$
Summing binomial coefficients using `Pascal's triangle' yields
\bea \label{pascal}
\pi^m_i=\left\{\begin{array}{ll}
{n\choose i}\left(\frac{1}{2}\right)^{n-1},&\ i\in L_{m,n},\\
0,&\ i\not\in L_{m,n}.
\end{array}\right.
\ena
We say that the random variable $C_{m,n}$ has the clubbed binomial distribution if $P(C_{m,n}=i)=\pi^m_i$ for $i \in L_{m,n}$. In words, the clubbed binomial distribution is formed by combining two adjacent cells of
the binomial.

It was observed in ~\cite{rrz} that the clubbed binomial distribution appeared to approximate the
lightbulb distribution $W_n$ exponentially well. Here we make that observation rigorous by supplying an exponentially decaying bound in total variation. First, recall that if $X$ and $Y$ are two random variables with distributions supported on $\mathbb{Z}$, then the total variation distance between the (laws of) $X$ and $Y$, denoted $d_{{\mbox{\scriptsize \rm TV}}}(X,Y)$, is given by
\bea \label{def:dtv}
d_{{\mbox{\scriptsize \rm TV}}}(X,Y)=\sup_{A \subset \mathbb{Z}}|P(X \in A)-P( Y \in A)|.
\ena

\begin{theorem}
\label{thm:illuminate}
Let $W_n$ be the total number of bulbs on at the terminal time in the lightbulb process of size $n$ and let $C_n=C_{m,n}$
where $m=0$ for  $n\,{\rm mod}\,4\in\{0,3\}$ and $m=1$ for $n\,{\rm mod}\,4\in\{1,2\}$. Then
\beas
d_{{\mbox{\scriptsize \rm TV}}}(W_n,C_n) \le 2.7314 \sqrt{n} e^{-(n+1)/3}.
\enas
\end{theorem}
In particular, the approximation error is less than 1\% for $n\ge 21$ and less than 0.1\% for $n\ge 28$.

A Berry-Esseen bound in the Kolmogorov metric of order $1/\sqrt{n}$ for the distance between the standardized value of $W_n$ and the unit normal was derived in ~\cite{gz}. The lighbulb chain was also studied in ~\cite{zl}, and served there as a basis for the exploration of the more general class of Markov chains of multinomial type. One feature of such chains is their easily obtainable spectral decomposition, which informed the analysis in \cite{gz}. In contrast, here we demonstrate the exponential bound in total variation using only simple properties of the lightbulb process.

After formalizing the framework for the lightbulb process in the next section, we prove Theorem \ref{thm:illuminate} by Stein's method. In particular, we develop a Stein operator ${\cal A}$ for the clubbed binomial distribution and obtain bounds on the solution $f$ of the associated Stein equation. The exponentially small distance between $W_n$ and the clubbed binomial $C_n$ can then be seen to be a consequence of the vanishing of the expectation of ${\cal A}f$ except on a set of exponentially small probability.

\section{The lightbulb process}

We now more formally describe the lightbulb process. With $n \in
\mathbb{N}$ fixed we will let ${\bf X}=\{X_{rk}: r=0,1,\ldots,n ,k=1,\ldots,n\}$ denote a collection of
Bernoulli variables. For $r \ge 1$ these `switch' or `toggle' variables
have the interpretation that
\beas X_{rk} &=& \left\{
\begin{array}{cc}
1 & \mbox{ if the status of bulb $k$ is changed at stage $r$,} \\
0 &\mbox{ otherwise.}
\end{array}
\right.
\enas
We take the initial state of the bulbs to be given deterministically by setting the switch variables $\{X_{0k},k=1,\ldots,n\}$ equal to zero, that is, all bulbs begin in the off position.
At stage $r$ for $r=1,\ldots,n$, $r$ of the $n$ bulbs are chosen uniformly to have their status changed, with
different stages mutually independent. Hence, with $e_1,\ldots,e_n \in \{0,1\}$,
the joint distribution of $X_{r1},\ldots,
X_{rn}$ is given by
\beas 
\lefteqn{P(X_{r1} = e_1, \cdots, X_{rn}=e_n)
=\left\{
\begin{array}{cc}
{n \choose r}^{-1} & \mbox{if $e_1+\cdots+e_n = r,$}  \\
0 & \mbox{ otherwise,}
\end{array}
\right.}\\
&& \nn \mbox{with the collections $\{X_{r1},\ldots,X_{rn}\}$
independent for $r=1,\ldots,n$.}
\enas
Clearly, at each stage $r$
the variables $(X_{r1},\cdots, X_{rn})$ are exchangeable.

For $r,i = 1, \ldots, n$, the
quantity $\left(\sum_{s=1}^r X_{si}\right) \mbox{ mod }2$ is the
indicator that bulb $i$ is on at time $r$ of the lightbulb process, so letting
\beas 
I_i = \left(\sum_{r=0}^n X_{ri}\right) \mbox{ mod }2
\qmq{and} W_n=\sum_{i=1}^n I_i,
\enas
the variable $I_i$ is the
indicator that bulb $i$ is on at the terminal time, and $W_n$ is the number of bulbs on at the terminal time.

The lightbulb process is a special case of a class of multivariate chains
studied in \cite{zl}, where randomly chosen subsets of $n$ individual particles evolve according
to the same marginal Markov chain. As shown in ~\cite{zl},
such chains admit explicit full spectral
decompositions, and in particular, the transition matrices for each stage of the lightbulb
process can be simultaneously diagonalized by a Hadamard matrix. These properties
were applied in \cite{rrz} for the calculation of the
moments needed to compute the mean and variance of $W_n$ and to develop
recursions for the exact distribution, and in \cite{gz} for a Berry-Esseen bound of
the standardized $W_n$ to the normal.

\section{Stein Operator}
In order to apply Stein's method, we first develop a Stein equation for the clubbed binomial distribution $C_{m,n}$  and then present bounds on its solution. With $\pi^m_x$ given by (\ref{pascal}), let
$\pi^m(A)=\sum_{x \in A} \pi^m_x$. Set $\alpha_x=(n-x)(n-1-x)$ and $\beta_x=x(x-1)$ for $x\in\{0,\dots,n\}$. One may easily directly verify the balance equation
\bea \label{eqn:balance}
\alpha_{x-2}\pi^m_{x-2} = \beta_x \pi^m_x \quad \mbox{for $x \in L_{m,n}$},
\ena
which gives the generator of the distribution
of $C_{m,n}$ as
\bea \label{generator}
{\cal A}f(x)=\alpha_x f(x+2)-\beta_x f(x),\mbox{ for } x\in L_{m,n}.
\ena

For $A \subset L_{m,n}$, we consider the Stein equation
\bea \label{Afr}
{\cal A}f_A(x)=1_{A}(x)-\pi^m(A),\ x\in L_{m,n}.
\ena
For a function $g$ with domain $A$ let $\|g\|$ denote $\sup_{x \in A} |g(x)|$.

\begin{lemma}
\label{lem:bound}
For $m\in\{0,1\}$ and $A=\{r\}$ with $r \in L_{m,n}$, the unique solution $f^m_r(x)$ of (\ref{Afr}) on $L_{m,n}$ satisfying the boundary condition $f^m_r(m)=0$ is given, for $m<x\le n, x \in L_{m,n}$, by
\bea \label{frsol}
f^m_r(x)= \left\{
\begin{array}{cl}
-\frac{\pi^m([0,x-2]\cap L_{m,n})\pi^m_r}{\beta_x \pi^m_x} & \mbox{for $m < x<r+2$} \\
\frac{\pi^m([x,n]\cap L_{m,n})\pi^m_r}{\beta_x \pi^m_x} &  \mbox{for $r+2 \le x \le n$.}
\end{array}
\right.
\ena
Furthermore, for all $A \subset L_{m,n}$, $f^m_A(x)=\sum_{r\in A}f^m_r(x)$ is a solution of (\ref{Afr}) and satisfies
\beas
\|f^m_A\| \le \frac{2.7314}{\sqrt{n}(n-1)} \quad \mbox{for $n \ge 1$.}
\enas
\end{lemma}
Lemma \ref{lem:bound} is proved in Section \ref{sec:bound}.

Applying Lemma \ref{lem:bound}, we now prove our main result.

{\em Proof of Theorem \ref{thm:illuminate}:} Fix $m\in\{0,1\}$ and $A \subset L_{m,n}$, and let $f:=f^m_A$ be the solution to (\ref{Afr}). Dropping subscripts,
let $W=\sum_{i=1}^n I_i$, where $I_i$ is the indicator that bulb $i$ is on at the terminal time. For $i,j \in \{1,\ldots,n\}$, now with slight abuse of notation, let $W_i=W-I_i$, and for $i \not = j$ set $W_{ij}=W-I_i-I_j$.
Then
\beas 
\lefteqn{E(n-W)(n-1-W)f(W+2)}\\
&=&E\sum_{i=1}^n (1-I_i)(n-1-W)f(W_i+2)\\
&=&E\sum_{i \not =j} (1-I_i)(1-I_j)f(W_{ij}+2),
\enas
and similarly,
\beas 
E W(W-1)f(W) = E\sum_{i=1}^n I_i W_i f(W_i+1) = E\sum_{i \not = j} I_i I_j f(W_i+1)= E\sum_{i \not = j} I_i I_j f(W_{ij}+2).
\enas

By Proposition~2 of \cite{rrz}, $P(W \in L_{m,n})=1$,
and hence (\ref{Afr}) holds upon replacing $x$ by $W$. Taking expectation and using the expression for the generator in (\ref{generator}), we obtain
\bea \label{gen-indicators}
P(W \in A) - \pi^m(A)
= E{\cal A}f(W)= E \sum_{i \not = j} \left( (1-I_i)(1-I_j)-I_iI_j \right) f(W_{ij}+2).
\ena
Recalling that $X_{rk}$ is the value of the switch variable at time $r$ for bulb $k$, let $A_{ij}$ be the event that the switch variables of the distinct bulbs $i$ and $j$ differ in at least one stage, that is, let
\bea \label{def:Aij}
A_{ij}= \bigcup_{r=1}^n \{X_{ri} \not = X_{rj}\}.
\ena
Now using (\ref{gen-indicators}) we obtain
\bea
&& \left| P(W \in A) - \pi^m(A) \right| = \left| E \sum_{i \not = j}  \left( (1-I_i)(1-I_j)-I_iI_j \right) f(W_{ij}+2) \right| \nn \\
&\le& \left| \sum_{i \not = j} E \left( (1-I_i)(1-I_j)-I_iI_j  \right) f(W_{ij}+2){\bf 1}_{A_{ij}} \right| \nn \\
&& + \left| \sum_{i \not = j} E \left( (1-I_i)(1-I_j)-I_iI_j \right)f(W_{ij}+2){\bf 1}_{A_{ij}^c}  \right|. \label{twoterms}
\ena
Note that $I_i, I_j \in \{0,1\}$ implies
$$
(1-I_i)(1-I_j){\bf 1}_{I_i \not = I_j}=0=I_iI_j {\bf 1}_{I_i \not = I_j},
$$
and hence for the first term in (\ref{twoterms}) we obtain
\bea
\lefteqn{\sum_{i \not = j} \left( (1-I_i)(1-I_j)-I_iI_j  \right) f(W_{ij}+2) {\bf 1}_{A_{ij}}}\nn \\
&=& \sum_{i \not = j} \left( (1-I_i)(1-I_j)-I_iI_j  \right) f(W_{ij}+2){\bf 1}_{A_{ij},I_i=I_j}. \label{eq:fterm}
\ena

For a given pair $i,j$, on the event $A_{ij}$ let $t$ be any index for which $X_{ti} \not = X_{tj}$, and let ${\bf X}^{ij}$ be the collection of switch variables given by
\beas
X_{rk}^{ij}= \left\{
\begin{array}{cl}
X_{rk} & r \not = t,\\
X_{tk} & r=t, k \not \in \{i,j\},\\
X_{ti} & r=t, k=j,\\
X_{tj} & r=t, k=i.
\end{array}
\right.
\enas
In other words, in stage $t$, the unequal switch variables $X_{ti}$ and $X_{tj}$ are interchanged, and all other
variables are left unchanged. Let $I_k^{ij}$ be the status of bulb $k$ at the terminal time when applying switch variables ${\bf X}^{ij}$, and similarly set
$W_{ij}^{ij}=\sum_{k \not \in \{i,j\}} I_k^{ij}$. Note that as the status of both bulbs $i$ and $j$ are toggled upon interchanging their stage $t$ switch variables, and all other variables are unaffected, we obtain
$$
I_i^{ij}=1-I_i, \quad I_j^{ij}=1-I_j \qmq{and} W_{ij}^{ij}=W_{ij}.
$$
In particular, $I_i=I_j$ if and only if $I_i^{ij}=I_j^{ij}$, and, with $A_{ij}^{ij}$ as in
(\ref{def:Aij}) with $X_{rk}^{ij}$ replacing $X_{rk}$, we have additionally that $A_{ij}^{ij}=A_{ij}$. Further, by exchangeability
 we have ${\cal L}({\bf X})={\cal L}({\bf X}^{ij})$. Therefore,
\beas
\lefteqn{E (1-I_i)(1-I_j)f(W_{ij}+2){\bf 1}_{A_{ij},I_i=I_j}}\\
&=&E (1-I_i^{ij})(1-I_j^{ij})f(W_{ij}^{ij}+2){\bf 1}_{A_{ij}^{ij},I_i^{ij}=I_j^{ij}}\\
&=&E I_i I_j f(W_{ij}+2){\bf 1}_{A_{ij},I_i=I_j},
\enas
showing, by (\ref{eq:fterm}), that the first term in (\ref{twoterms}) is zero. Therefore,
\beas
\lefteqn{|P(W \in A) - \pi^m(A) |}\\
&&\le \left| \sum_{i \not = j} E \left( (1-I_i)(1-I_j)-I_iI_j  \right) f(W_{ij}+2){\bf 1}_{A_{ij}^c} \right|
\le \|f\| \sum_{i \not = j}P(A_{ij}^c).
\enas
As $A_{ij}^c$ is the event that the switch variables of $i$ and $j$ are equal in every stage, recalling that
these variables are independent over stages we obtain
\beas
P(A_{ij}^c) &=& \prod_{r=1}^n \frac{r(r-1)+(n-r)(n-1-r)}{n(n-1)}\\
&=& \prod_{r=1}^n \left( 1 - \frac{2(nr-r^2)}{n(n-1)} \right) \\
&\le& e^{- \frac{2}{n(n-1)} \sum_{r=1}^n (nr-r^2)}
=e^{-(n+1)/3}.
\enas
Hence, by Lemma \ref{lem:bound},
\beas
\left|  P(W \in A) - \pi^m(A)  \right| \le \frac{2.7314}{\sqrt{n}(n-1)}n(n-1)e^{-(n+1)/3}= 2.7314 \sqrt{n} e^{-(n+1)/3}.
\enas
Taking supremum over $A$ and applying definition (\ref{def:dtv}) completes the proof.
\bbox

\section{Bounds on the Stein equation}
\label{sec:bound}
In this section we present the proof of Lemma \ref{lem:bound}.

\proof Let $m \in \{0,1\}$ be fixed. First, the equalities $f(m)=0$ and
\beas
f(x+2)=\frac{1_A(x)-\pi^m(A)+\beta_x f(x)}{\alpha_x} \quad \mbox{for $m<x \le n-2, x \in L_{m,n}$}
\enas
specify $f(x)$ on $L_{m,n}$ uniquely, hence the solution to (\ref{Afr}) satisfying the given boundary condition is unique.

Next, with $r \in L_{m,n}$, we verify that $f^m_r(x)$ given by (\ref{frsol}) solves (\ref{Afr}) with $A=\{r\}$;
that $f^m_r(m)=0$ is given.
For $m<x<r,x \in L_{m,n}$, applying the balance equation (\ref{eqn:balance}) to obtain the second equality,
we have
\beas
\lefteqn{\alpha_x f^m_r(x+2)-\beta_x f^m_r(x) }\\
&=& \alpha_x \left( -\frac{\pi^m([0,x]\cap L_{m,n})\pi^m_r}{\beta_{x+2} \pi^m_{x+2}} \right) - \beta_x \left( - \frac{\pi^m([0,x-2]\cap L_{m,n})\pi^m_r}{\beta_x \pi^m_x} \right)\\
&=& \alpha_x \left( -\frac{\pi^m([0,x]\cap L_{m,n})\pi^m_r}{\alpha_x \pi^m_x} \right) - \beta_x \left( - \frac{\pi^m([0,x-2]\cap L_{m,n})\pi^m_r}{\beta_x \pi^m_x} \right)\\
&=&-\pi^m_r.
\enas

If $x=r$ then
\beas
\lefteqn{\alpha_x f^m_r(x+2)-\beta_x f^m_r(x)}\\
&=& \alpha_r \left( \frac{\pi^m([r+2,n]\cap L_{m,n})\pi^m_r}{\beta_{r+2} \pi^m_{r+2}}\right)-\beta_r \left( \frac{-\pi^m([0,r-2]\cap L_{m,n})\pi^m_r}{\beta_r \pi^m_r}\right)\\
&=& \alpha_r \left(\frac{\pi^m([r+2,n]\cap L_{m,n})\pi^m_r}{\alpha_r \pi^m_r}\right)-\beta_r \left(\frac{-\pi^m([0,r-2]\cap L_{m,n})\pi^m_r}{\beta_r \pi^m_r}\right)\\
&=& \pi^m([r+2,n]\cap L_{m,n})+\pi^m([0,r-2]\cap L_{m,n})=1-\pi^m_r.
\enas

If $x>r$ then
\beas
\lefteqn{\alpha_x f^m_r(x+2)-\beta_x f^m_r(x)}\\
&=& \alpha_x \left(\frac{\pi^m([x+2,n]\cap L_{m,n})\pi^m_r}{\beta_{x+2} \pi^m_{x+2}}\right)-\beta_x \left(\frac{\pi^m([x,n]\cap L_{m,n}) \pi^m_r}{\beta_x \pi^m_x}\right)\\
&=& \alpha_x \left(\frac{\pi^m([x+2,n]\cap L_{m,n})\pi^m_r}{\alpha_x \pi^m_x}\right)-\beta_x \left(\frac{\pi^m([x,n]\cap L_{m,n}) \pi^m_r}{\beta_x \pi^m_x}\right)\\
&=&-\pi^m_r.
\enas
Hence $f^m_r(x)$ solves (\ref{Afr}).

Next, to consider the solution of (\ref{Afr}) more generally for $A \subset L_{m,n}$ and $x \in L_{m,n}$, letting
\beas
U_{m,x} = [0,x-2]\cap L_{m,n} \qmq{and} U_{m,x}^c = L_{m,n} \setminus U_{m,x},
\enas
we may write (\ref{frsol}) more compactly as
\beas
f^m_r(x)=\frac{1}{\beta_x \pi^m_x}\left( \pi^m(U_{m,x}^c) \pi^m(\{r\} \cap U_{m,x}) - \pi^m(U_{m,x}) \pi^m(\{r\} \cap U_{m,x}^c) \right).
\enas
By linearity, the solution of (\ref{Afr}) for $A \subset L_{m,n}$ is given by $f^m_A(m)=0$, and for $x>m,x \in L_{m,n}$, 
by
\beas 
f^m_A(x)=\frac{1}{\beta_x \pi^m_x}\left(\pi^m(U_{m,x}^c)\pi^m(A \cap U_{m,x})-\pi^m(U_{m,x})\pi^m(A \cap U_{m,x}^c) \right)
\enas
(cf \cite{bhj}, p. 7), and so, for all $x \in L_{m,n}$,
\beas
-\frac{1}{\beta_x \pi^m_x}\pi^m(U_{m,x})\pi^m(U_{m,x}^c) \le f_A^m(x) \le \frac{1}{\beta_x \pi^m_x}\pi^m(U_{m,x}^c) \pi^m(U_{m,x}),
\enas
or that
\bea
\left| f^m_A(x) \right| \le \frac{1}{\beta_x \pi^m_x} \pi^m(U_{m,x}) \pi^m(U_{m,x}^c).\label{smalln}
\ena

Since $f^m_A(m)=0$ and the upper bound of Lemma~\ref{lem:bound} reduces to $\infty$ if $0 \le n\le 1$, we only need to bound $f^m_A(x)$ for $n\ge 2$ and $x\ge 2$. Direct computation using (\ref{smalln}) gives $|f_A^0(2)|\le 1/4$ for $n=2$, $|f_A^0(2)|\le 1/8$ and $|f_A^1(3)|\le 1/8$ for $n=3$, $|f_A^0(2)|=|f_A^0(4)|\le 7/96$ and  $|f_A^1(3)|\le 1/12$ for $n=4$. Therefore, it remains to prove Lemma~\ref{lem:bound} for $n\ge 5$.

Noting that for $x\ge\frac n2+1$ we have $\beta_x\ge\left(\frac n2+1\right)\frac n2$, and for $x<\frac n2+1$ that $\alpha_{x-2}=(n-x+2)(n-x+1)>\left(\frac n2+1\right)\frac n2$,
using (\ref{eqn:balance}), we obtain from (\ref{smalln}) that
\bea \label{bnd:by-fAx}
\left| f^m_A(x) \right| \le \left\{\begin{array}{ll}
\frac{\pi^m(U_{m,x}) \pi^m(U_{m,x}^c)}{\beta_x \pi^m_x} \le \frac{1}{\left( \frac{n}{2}+1 \right)\frac
{n}{2}}\frac{\pi^m(U_{m,x}) \pi^m(U_{m,x}^c)}{\pi^m_x}&\mbox{ if }x\ge \frac n2+1,\\
\frac{\pi^m(U_{m,x}) \pi^m(U_{m,x}^c)}{\alpha_{x-2} \pi^m_{x-2}}\le \frac{1}{\left( \frac{n}{2}+1 \right)\frac
{n}{2}}\frac{\pi^m(U_{m,x}) \pi^m(U_{m,x}^c)}{\pi^m_{x-2}}&\mbox{ if }x< \frac n2+1.
\end{array}\right.
\ena

Clearly, for $i \ge x$,
\beas
\frac{\pi^m_i}{\pi^m_x}=\frac{{n \choose i}}{{n \choose x}} =
\left\{
\begin{array}{cl}
1 & \mbox{if $i=x$}\\
\frac{(n-x) \cdots (n-i+1)}{(x+1) \cdots i} & \mbox{if $i \ge x+2$.}
\end{array}
\right.
\enas
Hence, we can write, for $i \ge x+2$,
\bea \label{prodforigexp2}
\frac{\pi^m_i}{\pi^m_x}= \left( \frac{n-x}{x+1}\right)\left( \frac{n-x-1}{x+2}\right) \cdots \left( \frac{n-i+1}{i} \right)=\prod_{y=0}^{i-x-1}\frac{n-x-y}{x+1+y}.
\ena
Note that as $(n-x)/(x+1) \le 1$ for $x \ge n/2$, the terms in the product (\ref{prodforigexp2}) are decreasing. In particular,
\bea \label{star6}
\frac{\pi^m_i}{\pi^m_x} \le 1 \qmq{for $i \ge x$, and} \prod_{0 \le y \le
\lfloor \frac{\sqrt{n}}{2} \rfloor }\frac{n-x-y}{x+1+y} \le 1\, \mbox{ provided } x \ge \frac{n}{2}.
\ena

For $n$ even let $x_s=n/2$, and for $n$ odd let $x_s=(n-1)/2$ when $m=0$, and $x_s=(n+1)/2$ when $m=1$. Then, except for the case where $m=0$ and $x=(n+1)/2$, which we deal with separately, we have
$$
\pi^m(U_{m,x})\pi^m(U_{m,x}^c)=\pi^m(U_{m,2x_s-x+2})\pi^m(U_{m,2x_s-x+2}^c),
$$
and we may therefore assume $x \ge x_s+1$, and so $x \ge n/2+1$.

Since for $y \ge \sqrt{n}/2$, recalling $x \ge n/2+1$,  we have
\bea\label{star7}
\frac{n-x-y}{x+1+y} \le \frac{n-\left( \frac{n}{2}+1\right)-\frac{\sqrt{n}}{2}}{\left( 1+\frac{n}{2}\right)+1+\frac{\sqrt{n}}{2}}=\frac{\frac{n}{2}-\frac{\sqrt{n}}{2}-1}{\frac{n}{2}+2+\frac{\sqrt{n}}{2}}=1-\frac{\sqrt{n}+3}{\frac{n}{2}+2+\frac{\sqrt{n}}{2}},
\ena
applying (\ref{prodforigexp2}) and (\ref{star6}) we conclude that
\beas
\frac{\pi^m_i}{\pi^m_x} \le  \left(1-\frac{\sqrt{n}+3}{\frac{n}{2}+2+\frac{\sqrt{n}}{2}} \right)^{i-x-\lfloor \frac{\sqrt{n}}{2} \rfloor-1} \quad \mbox{for $i \ge x+\lfloor \frac{\sqrt{n}}{2} \rfloor+1$}.
\enas
Hence, applying (\ref{star6}) again, here to obtain the second inequality, we have
\bea
\lefteqn{\frac{1}{\left( \frac{n}{2}+1\right) \frac{n}{2}} \frac{\pi^m(U_{m,x}^c)}{\pi^m_x}}\nn \\
&\le&
\frac{1}{\left(\frac{n}{2}+1 \right) \frac{n}{2}} \left( \sum_{x \le i \le x + \lfloor\frac{\sqrt{n}}{2}\rfloor,i\in L_{m,n}} \frac{\pi^m_i}{\pi^m_x} + \sum_{i \ge x + \lfloor\frac{\sqrt{n}}{2}\rfloor+1,i\in L_{m,n}} \frac{\pi^m_i}{\pi^m_x} \right) \nn \\
&\le& \frac{1}{\left( \frac{n}{2} + 1\right) \frac{n}{2}} \left( \left( \frac{\sqrt{n}}{4}+1\right) + \sum_{j=0}^\infty \left( 1- \frac{\sqrt{n}+3}{\frac{n}{2}+2+\frac{\sqrt{n}}{2}}\right)^{2j} \right)\nn \\
&=& \frac{1}{\left(\frac{n}{2}+1 \right) \frac{n}{2}} \left(\frac{\sqrt{n}}{4}+1+\frac{1}{1-\left(1-\frac{\sqrt{n}+3}{\frac{n}{2}+2+\frac{\sqrt{n}}{2}}\right)^2} \right)\nn\\
&\le& \frac{2.7314}{\sqrt{n}(n-1)} \quad \mbox{for $n \ge 1$.}\label{g1n}
\ena
This final inequality is obtained by determining the maximum of the function
$$
g_1(n):=\frac{1}{\left(\frac{n}{2}+1 \right) \frac{n}{2}} \left(\frac{\sqrt{n}}{4}+1+\frac{1}{1-\left(1-\frac{\sqrt{n}+3}{\frac{n}{2}+2+\frac{\sqrt{n}}{2}}\right)^2} \right)\sqrt{n}(n-1)
$$
by noting $g_1(n)<1+\frac4{\sqrt{n}}+\frac{(n+4+\sqrt{n})^2}{(n+2)(n+3\sqrt{n})}<2.5$ for $n\ge 64$ and $\max_{1\le n\le 63}g_1(n)=g_1(9)=2.7313131\ldots$.

Lastly we handle the situation where $n$ is odd, $m=0$ and $x=(n+1)/2=:x_0 \in L_{0,n}$, in which case $n=3\,{\rm mod}\,4$.  In place of (\ref{star7}), we have, for $y \ge \sqrt{n}/2$,
\beas
\frac{n-x_0-y}{x_0+1+y} \le \frac{n-\left( \frac{n+1}{2}\right)-\frac{\sqrt{n}}{2}}{\left( \frac{n+1}{2}\right)+1+\frac{\sqrt{n}}{2}}=1-\frac{4+2\sqrt{n}}{n+3+\sqrt{n}}.
\enas
Since (\ref{star6}) is valid for all $x \ge n/2$, in view of (\ref{prodforigexp2}) we obtain the bound
\beas
\frac{\pi^m_i}{\pi^m_{x_0}} \le  \left(1-\frac{4+2\sqrt{n}}{n+3+\sqrt{n}} \right)^{i-x_0-\lfloor \frac{\sqrt{n}}{2} \rfloor-1} \quad \mbox{for $i \ge x_0+\lfloor \frac{\sqrt{n}}{2} \rfloor+1$}.
\enas
Using (\ref{star6}) again for the first inequality we have
\bea
\lefteqn{\frac{1}{\left( \frac{n}{2}+1\right) \frac{n}{2}} \frac{\pi^m(U_{m,x_0})\pi^m(U_{m,x_0}^c)}{\pi^m_{x_0}}}\nn \\
&=&
\frac{1}{2\left(\frac{n}{2}+1 \right) \frac{n}{2}} \left( \sum_{x_0 \le i \le x_0 + \lfloor\frac{\sqrt{n}}{2}\rfloor,i\in L_{m,n}} \frac{\pi^m_i}{\pi^m_{x_0}} + \sum_{i \ge x_0 + \lfloor\frac{\sqrt{n}}{2}\rfloor+1,i\in L_{m,n}} \frac{\pi^m_i}{\pi^m_{x_0}} \right) \nn \\
&\le& \frac{1}{\left( \frac{n}{2} + 1\right) n} \left( \left( \frac{\sqrt{n}}{4}+1\right) + \sum_{j=0}^\infty  \left(1-\frac{4+2\sqrt{n}}{n+3+\sqrt{n}} \right)^j \right) \nn\\
&=& \frac{1}{\left(\frac{n}{2}+1 \right) n} \left(\frac{\sqrt{n}}{4}+1+\frac{n+3+\sqrt{n}}{4+2\sqrt{n}} \right)\nn\\
&\le& \frac{1.638496535}{\sqrt{n}(n-1)} \quad \mbox{for $n \ge 1$,}\label{g2n}
\ena
where the last inequality is from bounding the function
$$g_2(n):=\frac{1}{\left(\frac{n}{2}+1 \right) n} \left(\frac{\sqrt{n}}{4}+1+\frac{n+3+\sqrt{n}}{4+2\sqrt{n}}\right)(n-1)\sqrt{n},$$
with $g_2(n)\le \frac{1}{2}+\frac{2}{\sqrt{n}}+\frac{n+3+\sqrt{n}}{n+2\sqrt{n}}\le 1.6$ for $n\ge 400$ and $\max_{1\le n\le 399}g_2(n)=g_2(23)=1.638496535$.

The result now follows from combining the estimates (\ref{g1n}), (\ref{g2n}) and (\ref{bnd:by-fAx}).
\bbox

We remark that a direct argument using Stirling's formula for the case $x=\lfloor n/2\rfloor$ shows that the best order that can be achieved for the estimate of $f^m_A$ is $O(n^{-3/2})$.



{\bf Acknowledgement:} The authors would like to thank the organizers of the conference held at the National University of Singapore in honor of Louis Chen's birthday for the opportunity to collaborate on the present work.

\end{document}